\documentclass{article}
\usepackage{amssymb}
\usepackage{amsmath}
\usepackage[all]{xy}
\title{\Large \bf Syzygy modules for quasi $k$-Gorenstein rings
}
\author{Zhaoyong Huang
\thanks{{\small \it E-mail address: huangzy@nju.edu.cn}} \\
{\small \it Department of Mathematics, Nanjing University,}\\
{\small \it Nanjing 210093, People's Republic of China}}
\usepackage{amssymb}
\date{}
\begin{document}
\baselineskip=18pt \maketitle

\begin{abstract}
Let $\Lambda$ be a quasi $k$-Gorenstein ring. For each $d$th
syzygy module $M$ in mod $\Lambda$ (where $0 \leq d \leq k-1$), we
obtain an exact sequence $0 \to B \to M \bigoplus P \to C \to 0$
in mod $\Lambda$ with the properties that it is dual exact, $P$ is
projective, $C$ is a $(d+1)$st syzygy module, $B$ is a $d$th
syzygy of Ext$_{\Lambda}^{d+1}(D(M), \Lambda)$ and the right
projective dimension of $B^*$ is less than or equal to $d-1$. We
then give some applications of such an exact sequence as follows.
(1) We obtain a chain of epimorphisms concerning $M$, and by
dualizing it we then get the spherical filtration of Auslander and
Bridger for $M^*$. (2) We get Auslander and Bridger's
Approximation Theorem for each reflexive module in mod $\Lambda
^{op}$. (3) We show that for any $0 \leq d \leq k-1$ each $d$th
syzygy module in mod $\Lambda$ has an Evans-Griffith
representation. As an immediate consequence of (3), we have that,
if $\Lambda$ is a commutative noetherian ring with finite
self-injective dimension, then for any non-negative integer $d$,
each $d$th syzygy module in mod $\Lambda$ has an Evans-Griffith
representation, which generalizes an Evans and Griffith's result
to much more general setting.
\end{abstract}

\vspace{0.5cm}

\centerline{\bf 1. Introduction}

\vspace{0.25cm}

Let $\Lambda$ be a left and right noetherian ring and
mod $\Lambda$ the category of finitely generated left
$\Lambda$-modules.

It is well known that $\Lambda$ possesses rather interesting
properties when it satisfies the condition
that gradeExt$_{\Lambda}^i(N, \Lambda)\geq i$ for any
$N\in$mod $\Lambda ^{op}$ and $1 \leq i \leq k$ (where
$k$ is a positive integer). Assume that $\Lambda$ satisfies
this grade condition. For any $T$ in mod $\Lambda ^{op}$, Auslander and
Bridger in [3] Spherical Filtration Theorem 2.37 produced a
projective module $Q$ and a filtration $T\bigoplus Q=T_{0}\supseteq
T_{1} \supseteq \cdots \supseteq T_{k}$
in mod $\Lambda ^{op}$ such that each
$T_{i}/T_{i+1}$ is ``spherical" in the sense that the
cohomological Ext$_{\Lambda}^{j}(T_{i}/T_{i+1}, \Lambda)\neq 0$
only if $j=0$ or $j=i$. They also showed in [3] that
under this grade condition the following statements
are true: (1) The
full subcategory of mod $\Lambda ^{op}$ consisting of the modules
with projective dimension less than or equal to $k$ is covariantly
finite (see [4] for the definition of covariantly
finite); (2) A $d$th syzygy module in mod $\Lambda ^{op}$ is
$d$-torsionfree for any $1 \leq d \leq k$.
We remark that the second statement
doesn't hold in general although the converse is always true.

Under the above grade condition, Auslander and Reiten in [5]
proved that the right flat dimension of the $i$th term in a
minimal injective resolution of $\Lambda$ as a right
$\Lambda$-module is less than or equal to $i$ for any $1\leq i\leq
k$; recently, Hoshino and Nishida in [10] proved that the converse
also holds. We call a ring quasi $k$-Gorenstein provided it
satisfies one of these equivalent conditions. A ring is called
quasi Auslander ring if it is quasi $k$-Gorenstein for all $k$.

Recall that $\Lambda$ is called a $k$-Gorenstein ring if the right
flat dimension of the $i$th term in a minimal injective resolution
of $\Lambda$ as a right $\Lambda$-module is less than or equal to
$i-1$ for any $1\leq i\leq k$. This notion was introduced by
Auslander in [9]. Iwanaga and Sato in [11] called $\Lambda$ an
Auslander ring if it is $k$-Gorenstein for all $k$. In [9] Theorem
3.7 Auslander showed that the notion of $k$-Gorenstein rings (and
hence that of Auslander rings) is left-right symmetric and that
$\Lambda$ is $k$-Gorenstein if and only if the grade of any
submodule of Ext$_{\Lambda}^i(N, \Lambda)$ for any $N\in$mod
$\Lambda ^{op}$ and $1 \leq i \leq k$ is greater than or equal to
$i$.  However, as already pointed out in [5], the notion of quasi
$k$-Gorenstein rings (and hence that of quasi Auslander rings) is
not left-right symmetric.

Notice that Bass showed in [7] that a commutative noetherian ring
$\Lambda$ has finite self-injective dimension if and only if
gradeExt$_{\Lambda}^i(N, \Lambda)\geq i$ for any $N\in$mod
$\Lambda$ and $i \geq 1$. So the notion of Auslander rings is in
fact a generalization of that of commutative noetherian rings with
finite left and right self-injective dimensions.

The discussion in this paper is based on the results mentioned
above. For a quasi $k$-Gorenstein ring $\Lambda$ and each $d$th
syzygy module in mod $\Lambda$ (where $0 \leq d \leq k-1$) we
obtain here an exact sequence with ``nice" properties as follows.

\vspace{0.2cm}

{\bf Theorem} {\it Let} $\Lambda$ {\it be a quasi} $k$-{\it
Gorenstein ring and} $M$ {\it a} $d$th {\it syzygy module in} mod
$\Lambda$ ({\it where} $d$ {\it is an integer with} $0 \leq d \leq
k-1$). {\it Then there is a projective module} $P$ {\it in} mod
$\Lambda$ {\it such that the} $d$th {\it syzygy} $B$ {\it of}
Ext$_{\Lambda}^{d+1}(D(M), \Lambda)$ ({\it see Section 2 for the
definition of} $D(M)$) {\it is a submodule of} $M \bigoplus P$
{\it and such that the exact sequence}
$$0 \to B \to M \bigoplus P
\to C \to 0$$
{\it has the following properties}:

(1) $C$ {\it is a} $(d+1)$st {\it syzygy module}.

(2) r.pd$_{\Lambda}(B^{*})\leq d-1$.

(3) {\it The sequence} $0 \to B \to M\bigoplus P \to C
\to 0$ {\it is dual exact, that is, the induced sequence}
$0 \to C^* \to M^*\bigoplus P^* \to B^* \to 0$
{\it is exact}.

\vspace{0.2cm}

The above theorem is the main result in this paper, we will prove
it in Section 3. To prove it, we collect some preliminary results
in Section 2. In Section 4 we give some applications of our main
theorem. For example, as an application of the theorem, we obtain
a chain of epimorphisms concerning a module $M$ in mod $\Lambda$,
by dualizing it we then get the spherical filtration of Auslander
and Bridger for $M^*$; and furthermore we get Auslander-Bridger's
Approximation Theorem for each reflexive module in mod $\Lambda
^{op}$.

Evans and Griffith in [7] Theorem 2.1 showed that if $\Lambda$ is
a commutative noetherian local ring with finite global dimension
and contains a field then each non-free $d$th syzygy of rank $d$
has an Evans-Griffith representation. As another application of
Theorem above, we show that for a quasi $k$-Gorenstein ring
$\Lambda$ and any $0 \leq d \leq k-1$ each $d$th syzygy module in
mod $\Lambda$ has an Evans-Griffith representation; especially, we
have that, if $\Lambda$ is a commutative noetherian ring with
finite self-injective dimension, then for any non-negative integer
$d$, each $d$th syzygy module in mod $\Lambda$ has an
Evans-Griffith representation, which generalizes Evans and
Griffith's result in [7] Theorem 2.1 to much more general setting.

\vspace{0.5cm}

\centerline{\bf 2. Preliminaries}

\vspace{0.25cm}

In this section, we give some definitions in our terminology and
collect some facts which are used in this paper.

Throughout this paper, $\Lambda$ is a left and right noetherian
ring, mod $\Lambda$ is the category of finitely generated left
$\Lambda$-modules and $\Omega ^k({\rm mod}\ \Lambda)$ is the
full subcategory of mod $\Lambda$ consisting of $k$th syzygy modules.
Let $M$ be a module in mod $\Lambda$ (resp. mod $\Lambda ^{op}$).
We use l.pd$_{\Lambda}(M)$ (resp. r.pd$_{\Lambda}(M)$) to denote
the left (resp. right) projective dimension of $M$.
We use $\sigma _{M}: M \to M^{**}$, defined by
$\sigma _{M}(x)(f)=f(x)$ for any $x\in M$ and $f\in M^*$,
to denote the canonical evaluation homomorphism. $M$ is called
torsionless if $\sigma _{M}$ is a monomorphism; and $M$ is called
reflexive if $\sigma _{M}$ is an isomorphism.
For a non-negative integer $i$,
we denote grade$M\geq i$ if Ext$_{\Lambda}^j(M, \Lambda)=0$ for any
$0 \leq j <i$.

Let $M$ be in mod $\Lambda$ (resp. mod $\Lambda ^{op}$) and
$$P_{1}\to P_{0} \to M \to 0$$ a projective resolution of $M$ in mod
$\Lambda$ (resp. mod $\Lambda ^{op}$). Then we have an exact
sequence
$$0 \to M^{*} \to P_{0}^{*} \to P_{1}^{*} \to D(M) \to 0$$ in
mod $\Lambda ^{op}$ (resp. mod $\Lambda$), where $D(M)={\rm
Coker}(P_{0}^{*} \to P_{1}^{*})$. The following lemma is due to Auslander.

\vspace{0.2cm}

{\bf Lemma 2.1} ([2] Proposition 6.3) {\it Let} $M$ {\it and}
$D(M)$ {\it be as above. Then we have the following exact sequence}:
$$0 \to {\rm Ext}_{\Lambda}^1(D(M), \Lambda) \to M
\buildrel {\sigma _M} \over \longrightarrow M^{**} \to
{\rm Ext}_{\Lambda}^2(D(M), \Lambda) \to 0.$$

\vspace{0.2cm}

It is clear that Ext$_{\Lambda}^{i}(D(M),
\Lambda)\cong$Ext$_{\Lambda}^{i-2}(M^{*}, \Lambda)$ for any $i
\geq 3$. On the other hand, Ext$_{\Lambda}^{1}(D(M), \Lambda)$
\linebreak
$\cong
{\rm Ker}\sigma _{M}$ and Ext$_{\Lambda}^{2}(D(M), \Lambda)\cong {\rm
Coker}\sigma _{M}$ by Lemma 2.1. So we get that, although
$D(M)$ depends on the choice of the projective resolution of $M$, each of
Ext$_{\Lambda}^{i}(D(M), \Lambda)$ (for any $i\geq 1$) is independent
of the choice of the projective resolution of $M$ and hence is identical
up to isomorphisms.

Recall that $M$ is called $k$-torsionfree if
Ext$_{\Lambda}^i(D(M), \Lambda)=0$ for any $1 \leq i \leq k$
(see [3]). By Lemma 2.1, we have that $M$ is 1-torsionfree
(resp. 2-torsionfree) if and only if it is torsionless
(resp. reflexive). We use $\mathcal{T}^k({\rm mod}\ \Lambda)$
to denote the full subcategory of mod $\Lambda$ consisting of
$k$-torsionfree modules. It is easy to see that
$\mathcal{T}^k({\rm mod}\ \Lambda)\subseteq
\Omega ^k({\rm mod}\ \Lambda)$. Furthermore, we have the following
useful result, which gives some equivalent conditions of
$\mathcal{T}^i({\rm mod}\ \Lambda)=
\Omega ^i({\rm mod}\ \Lambda)$ for any $1 \leq i \leq k$.

\vspace{0.2cm}

{\bf Lemma 2.2} {\it For a positive integer} $k$, {\it the
following statements are equivalent}.

(1) gradeExt$_{\Lambda}^{i+1}(M, \Lambda)\geq i$ {\it for any}
$M\in$ mod $\Lambda$ {\it and} $1\leq i \leq k-1$.

(2) $\Omega ^i({\rm mod}\ \Lambda)=\mathcal{T}^i({\rm mod}\
\Lambda)$ {\it for any} $1\leq i \leq k$.

(3) gradeExt$_{\Lambda}^{i+1}(N, \Lambda)\geq i$ {\it for any}
$N\in$ mod $\Lambda ^{op}$ {\it and} $1\leq i \leq k-1$.

(4) $\Omega ^i({\rm mod}\ \Lambda ^{op})=\mathcal{T}^i({\rm mod}\
\Lambda ^{op})$ {\it for any} $1\leq i \leq k$.

\vspace{0.2cm}

{\it Proof.} The equivalence of (1) and (2) is proved in
[3] Proposition 2.26 (or see [5] Proposition 1.6). The
other implications are proved in [9] Theorem 2.4. $\blacksquare$

\vspace{0.2cm}

{\bf Corollary 2.3} {\it If} $\Lambda$ {\it is a quasi} $k$-{\it
Gorenstein ring, then} $\Omega ^i({\rm mod}\
\Lambda)=\mathcal{T}^i({\rm mod}\ \Lambda)$ {\it and} $\Omega
^i({\rm mod}\ \Lambda ^{op})= \mathcal{T}^i({\rm mod}\ \Lambda
^{op})$ {\it for any} $1\leq i \leq k+1$.

\vspace{0.2cm}

{\it Proof.} By [5] Proposition 4.2 and Theorem 1.7, we have that
gradeExt$_{\Lambda}^{i+1}(M, \Lambda)\geq i$ for any
$M\in$ mod $\Lambda$ and $1\leq i \leq k$. Now our conclusion follows
from Lemma 2.2. $\blacksquare$

\vspace{0.5cm}

\centerline{\bf 3. The Proof of The Theorem}

\vspace{0.25cm}

In this section, we will prove the theorem mentioned in Introduction. We
proceed in several steps.

\vspace{0.2cm}

{\it Proof of Theorem.} \underline{{\it The case} $d=0$}:
Put $C={\rm Im}\sigma _{M}$ and
$B=$Ext$_{\Lambda}^{1}(D(M), \Lambda)$.
Then we have an exact sequence in mod $\Lambda$:
$$0 \to B
\to M \to C \to 0.$$ Since $C$
is a submodule of $M^{**}$, $C$ is torsionless and $C\in \Omega
^{1}({\rm mod}$ $\Lambda)$. On the other hand,
gradeExt$_{\Lambda}^{1}(D(M), \Lambda) \geq 1$ by assumption, that is,
$B^*=0$, so
the obtained exact sequence is desired.

\underline{{\it The case} $d=1$}: Assume that
$$0 \to B \buildrel
{f} \over \longrightarrow P \buildrel {g} \over \longrightarrow
{\rm Ext}_{\Lambda}^{2}(D(M), \Lambda) \to 0$$
is an exact sequence in mod $\Lambda$ with
$P$ projective.
Consider the following pull-back diagram with the middle row
splitting:

$$\xymatrix{& & 0 \ar[d] & 0 \ar[d]& &\\
& & B \ar@{=}[r] \ar[d] & B \ar[d]^{f}& &\\
0 \ar[r] & M \ar@{=}[d] \ar[r] & M\bigoplus P
\ar[d] \ar[r] & P \ar[d]^{g} \ar[r] & 0\\
0 \ar[r] & M \ar[r]^{\sigma _{M}} &
M^{**} \ar[r] \ar[d] & {\rm Ext}^{2}_{\Lambda}(D(M),
\Lambda) \ar[d] \ar[r] & 0 &\\
& & 0 & 0 & & }
$$

Because $M^{**}$ is a dual, $M^{**}$ is a
second syzygy. Since gradeExt$_{\Lambda}^{2}(D(M), \Lambda)\geq 2$ by
assumption, $P^{*} \buildrel {f^{*}} \over \longrightarrow
B^{*}$ is an isomorphism and $B^{*}$ is projective. We know from
[1] Proposition 20.14 that $\sigma _{M}^{*}$ is
epic, so, by applying the functor Hom$_{\Lambda}(-, \Lambda)$ to
the above diagram, we get the following commutative diagram with
exact columns and rows:

$$\xymatrix{& 0 \ar[d] & 0 \ar[d]& & &\\
0 \ar[r] & [{\rm Ext}^{2}_{\Lambda}(D(M), \Lambda)]^{*}(=0) \ar[r]
\ar[d]^{g^{*}} & M^{***} \ar[d] \ar[r]^{\sigma
^{*}_{M}} & M^{*}
\ar@{=}[d] \ar[r] & 0\\
0 \ar[r] & P^{*} \ar[r] \ar[d]^{f^{*}} & P^{*}\bigoplus
M^{*} \ar[r] \ar[d] & M^{*} \ar[r] & 0\\
& B^{*} \ar[d] \ar@{=}[r] & B^{*} & & &\\
& 0 & & & & }
$$
It is easy to see that Coker$(M^{***}\to P^{*}\bigoplus
M^{*})=B^{*}$. Then the middle column in the former
diagram:
$$0\to B \to M\bigoplus P \to M^{**} \to 0$$
is desired.

\underline{{\it The case} $d\geq 2$}: Let
$$\cdots \to P_{i} \to \cdots \to
P_{1} \to P_{0} \to M^{*} \to 0
\eqno{(1)}$$ be a projective
resolution of $M^{*}$ in mod $\Lambda ^{op}$.
If $d=2$, then
$M$ is reflexive by Lemma 2.2. So we have an
exact sequence:

$$0\to M(\cong M^{**})
\to P_{0}^{*} \to P_{1}^{*} \to N \to 0,$$ where $N={\rm
Coker}(P_{0}^{*} \to P_{1}^{*})$. Now suppose $d\geq 3$. Since
$M$ is a $d$th syzygy module, It follows from
Lemma 2.2 that $M$ is $d$-torsionfree and so
Ext$_{\Lambda}^i(M^{*}, \Lambda)=0$ for any $1 \leq i \leq d-2$.
From this fact and the exact sequence (1)
we yield the following exact sequence:

$$0\to M(\cong M^{**})\to
P_{0}^{*} \to P_{1}^{*} \to \cdots \to P_{d-2}^{*} \to
P_{d-1}^{*}\to N \to 0
\eqno{(2)}$$ where $N={\rm
Coker}(P_{d-2}^{*} \to P_{d-1}^{*})$. So for any $d\geq 2$ we have
an exact sequence of the form (2).

By Lemma 2.1 we get easily the following exact sequence:

$$0 \to {\rm Ext}_{\Lambda}^{d+1}(D(M), \Lambda) \to N
\buildrel {\sigma _{N}} \over \longrightarrow  N^{**} \to {\rm
Ext}_{\Lambda}^{d+2}(D(M), \Lambda) \to 0.$$

Write $K={\rm Ext}_{\Lambda}^{d+1}(D(M), \Lambda)$ and $Y={\rm
Im}\sigma _{N}$. Let $\mathbb{U}^{\bullet}$ and
$\mathbb{V}^{\bullet}$ be the projective resolutions of $K$ and
$Y$, respectively. Then there is a projective module $P'$ in mod
$\Lambda$ such that we have the following commutative diagram with
exact columns and rows:

$$\xymatrix{& 0 \ar[d] & 0 \ar[d] & 0 \ar[d] &\\
0 \ar[r] & B \ar[d] \ar[r] & \underline{M}\bigoplus P'
\ar[d] \ar[r] & C' \ar[d] \ar[r]& 0\\
0 \ar[r] & U_{d-1} \ar[d] \ar[r] & U_{d-1}\bigoplus V_{d-1}
\ar[d] \ar[r] & V_{d-1} \ar[d] \ar[r]& 0\\
& \vdots \ar[d] & \vdots \ar[d] & \vdots \ar[d] &\\
 0 \ar[r] & U_{0} \ar[d]
\ar[r] & U_{0}\bigoplus V_{0}
\ar[d] \ar[r] & V_{0} \ar[d] \ar[r]& 0\\
0 \ar[r] & K \ar[d] \ar[r] & N
\ar[d] \ar[r] & Y \ar[d] \ar[r]& 0\\
& 0 & 0 & 0 & }
$$ where
$\underline{M}$ is the greatest direct summand
of $M$ without projective summands.

Since $Y$ is a submodule of $N^{**}$, $Y$ is
torsionless and hence it is in $\Omega ^{1}($mod $\Lambda)$.
So $C'$ is in $\Omega ^{d+1}($mod $\Lambda)$. On the other hand,
grade$K$=gradeExt$_{\Lambda}^{d+1}(D(M), \Lambda)\geq
d+1$ by assumption, so we get an exact sequence $0\to U_{0}^{*} \to
\cdots \to U_{d-1}^{*} \to B^{*} \to 0$ and hence
r.pd$_{\Lambda}(B^{*})\leq d-1$. It is trivial that every
homomorphism $f:B\to \Lambda$ may extends to a homomorphism $g:
U_{d-1}\to \Lambda$, so $f$ may extends to a
homomorphism $h:\underline{M}\bigoplus P' \to \Lambda$ and hence
the sequence $0 \to {C'}^{*} \to \underline{M}^{*}\bigoplus {P'}^{*}
\to B^{*} \to 0$ is exact.

Put $M=\underline{M}\bigoplus P''$, $P=P'\bigoplus P''$ and
$C=C'\bigoplus P''$. Then from the exact sequence
$0 \to B \to \underline{M}\bigoplus P' \to C' \to 0$ we yield
the following exact sequence:
$$0 \to B \to M\bigoplus P \to C \to 0,$$
which is desired. $\blacksquare$

\vspace{0.5cm}

\centerline{\bf 4. Applications}

\vspace{0.25cm}

In this section we will give some applications of the main theorem.
We first have the following result.

\vspace{0.2cm}

{\bf Corollary 4.1} {\it If} $\Lambda$ {\it is a quasi Auslander
ring, then for any non-negative integer} $d$ {\it and } $M$ {\it
in} $\Omega ^{d}({\rm mod}\ \Lambda)$, {\it there is a projective
module} $P$ {\it in} mod $\Lambda$ {\it such that the} $d$th {\it
syzygy} $B$ {\it of} Ext$_{\Lambda}^{d+1}(D(M), \Lambda)$ {\it is
a submodule of} $M\bigoplus P$ {\it and such that the exact
sequence} $0 \to B \to M\bigoplus P \to C \to 0$ {\it has the
following properties}:

(1) $C \in \Omega ^{d+1}({\rm mod}\ \Lambda)$.

(2) r.pd$_{\Lambda}(B^{*})\leq d-1$.

(3) {\it The sequence} $0 \to B \to M\bigoplus P \to C
\to 0$ {\it is dual exact.}

\vspace{0.2cm}

{\bf Proposition 4.2} ({\it The duality of spherical filtration})
{\it Let} $\Lambda$ {\it be a quasi} $k$-{\it Gorenstein ring.
Then, for each} $M$ {\it in} mod $\Lambda$, {\it there is a
projective module} $P$ {\it in} mod $\Lambda$ {\it and a chain of
epimorphisms:}
$$M \bigoplus P =M_0 \twoheadrightarrow M_1 \twoheadrightarrow
\cdots \twoheadrightarrow M_{k-1}\twoheadrightarrow M_k,$$
{\it such that}

(1) $B_d=$Ker$(M_d \to M_{d+1})$ {\it is a} $d$th {\it syzygy of}
Ext$_{\Lambda}^{d+1}(D(M), \Lambda)$
({\it or equivalently},
$B_d$ {\it is a} $d$th {\it syzygy of} Ext$_{\Lambda}^{d-1}(M^*, \Lambda)$
{\it if} $d \geq 2$) {\it for any} $0 \leq d \leq k-1$.

(2) $M_d \in \Omega ^d({\rm mod}\ \Lambda)$
{\it for any} $0 \leq d \leq k$.

(3) r.pd$_{\Lambda}(B_d^*)\leq d-1$ {\it for any} $0 \leq d \leq k-1$.

(4) {\it Each exact sequence} $0 \to B_d \to M_d \to M_{d+1} \to 0$
{\it is dual exact for any} $0 \leq d \leq k-1$.

\vspace{0.2cm}

{\it Proof.} We proceed by employing induction with successive applications
of Theorem in Introduction.

First, by Theorem and its proof, we have
an exact sequence in mod $\Lambda$:
$$0 \to B_0 \to M \to C_1 \to 0
\eqno{(3)}$$
with the properties that it is dual exact,
$B_0={\rm Ext}_{\Lambda}^1(D(M), \Lambda)$ and $C_1={\rm Im}\sigma _M$.
Notice that $B_0^*=0$ and
${\rm Im}\sigma _M$ is torsionless, then by Lemma 2.1 we have
that ${\rm Ext}_{\Lambda}^2(D(C_1), \Lambda)=
{\rm Ext}_{\Lambda}^2(D({\rm Im}\sigma _M), \Lambda) \cong
({\rm Im}\sigma _M)^{**}/{\rm Im}\sigma _M \cong
M^{**}/{\rm Im}\sigma _M \cong{\rm Ext}_{\Lambda}^2(D(M), \Lambda)$.

Next, by Theorem and its proof, we have
an exact sequence in mod $\Lambda$:
$$0 \to B_1 \to C_1 \bigoplus P_1 \to C_2 \to 0$$
with the properties that it is dual exact, $C_2=C_1^{**}$, $P_1$ is
projective, $B_1$ is a first syzygy of
Ext$_{\Lambda}^2(D(C_1), \Lambda)\cong$(Ext$_{\Lambda}^2(D(M), \Lambda)$),
$B_1^*$ is projective and $C_2 \in \Omega ^2({\rm mod}\ \Lambda)$.
Then we have that
$C_2^* \bigoplus B_1^* \cong C_1^* \bigoplus P_1^* \cong
M^* \bigoplus P_1^*$ and
Ext$_{\Lambda}^i(C_2^*, \Lambda)
\cong$Ext$_{\Lambda}^i(M^*, \Lambda)$ for any $i \geq 1$.

Now suppose that $k \geq 3$ and for any $0 \leq d \leq k-2$ there is
an exact sequence in mod $\Lambda$:
$$0 \to B_d \to C_d \bigoplus P_d \to C_{d+1} \to 0$$
with the properties that it is dual exact, $P_d$ is
projective, $B_d$ is a $d$th syzygy of
Ext$_{\Lambda}^{d+1}(D(C_d), \Lambda)$,
r.pd$_{\Lambda}B_d^* \leq d-1$, $C_0=M$ and
$C_{d+1} \in \Omega ^{d+1}({\rm mod}\ \Lambda)$.
Then we have that
Ext$_{\Lambda}^{k-2}(C_{k-1}^*, \Lambda)\cong$
Ext$_{\Lambda}^{k-2}(C_{k-2}^*, \Lambda)\cong \cdots \cong$
Ext$_{\Lambda}^{k-2}(C_2^*, \Lambda)\cong$
Ext$_{\Lambda}^{k-2}(M^*, \Lambda)$.

By Theorem, there is a projective
module $P_{k-1}$ and an exact sequence in mod $\Lambda$:
$$0 \to B_{k-1} \to C_{k-1} \bigoplus P_{k-1} \to C_k \to 0,$$
such that

(1) $B_{k-1}$ is a $(k-1)$st syzygy of Ext$_{\Lambda}^k
(D(C_{k-1}), \Lambda)$; or equivalently,
$B_{k-1}$ is a $(k-1)$st syzygy of
\linebreak
Ext$_{\Lambda}^{k-2}
(C_{k-1}^*, \Lambda)(\cong$Ext$_{\Lambda}^{k-2}
(M^*, \Lambda)\cong$Ext$_{\Lambda}^k
(D(M), \Lambda))$.

(2) $C_k \in \Omega ^k({\rm mod}\ \Lambda)$.

(3) r.pd$_{\Lambda}(B_{k-1}^*)\leq k-2$.

(4) The induced sequence $0 \to C_k^*
\to C_{k-1}^*\bigoplus P_{k-1}^*
\to B_{k-1}^* \to 0$ is exact.

Put $P=\bigoplus _{i=1}^{k-1}P_i$, $M_0=M\bigoplus P$,
$M_i=C_i\bigoplus (\bigoplus _{j=i}^{k-1} P_j)$
for any $1 \leq i \leq k-1$ and $M_k=C_k$. Then we get our conclusion.
$\blacksquare$

\vspace{0.2cm}

Let $T\in$mod $\Lambda ^{op}$. We remark that if one takes a chain
of epimorphisms of $T^*$ as in Proposition 4.2 and dualizes it,
one then obtains the spherical filtration of Auslander and Bridger
for $T^{**}$. Notice that a module $T$ in mod $\Lambda ^{op}$ is
reflexive if and only if it is isomorphic to $T^{**}$. So we in
fact obtain the spherical filtration of Auslander and Bridger for
each reflexive module in mod $\Lambda ^{op}$, and thus we may
regard Proposition 4.2 as a duality of the spherical filtration of
Auslander and Bridger.

As a corollary of Proposition 4.2, we get Auslander-Bridger's
Approximation Theorem (see [3] Theorem 2.41)
for $T^{**}$ (or for $T$ if
$T$ is reflexive) as follows.

\vspace{0.2cm}

{\bf Corollary 4.3} {\it Let} $\Lambda$ {\it be a quasi} $k$-{\it
Gorenstein ring. Then, for any} $T\in$mod $\Lambda ^{op}$, {\it
there are a projective module} $P$ {\it and an exact sequence in}
mod $\Lambda ^{op}$:
$$0 \to X \to T^{**}\bigoplus P \to Y \to 0$$
{\it such that}

(1) {\it It is dual exact}.

(2) r.pd$_{\Lambda}Y \leq k-2$.

(3) {\it The homomorphism} $T^{**}\bigoplus P \to Y$
{\it induces isomorphisms} ${\rm Ext}_{\Lambda}^i(Y, \Lambda)
\buildrel {\cong} \over \to
{\rm Ext}_{\Lambda}^i(T^{**}, \Lambda)$ for any $1 \leq i \leq k-2$.

(4) {\it If} $\underline{T^{**}} \to \underline{H}$ {\it is a homomorphism with}
r.pd$_{\Lambda}\underline{H} \leq k-2$, {\it then the above exact sequence
induces an isomorphism} ${\rm Hom}_{\Lambda}(\underline{Y}, \underline{H})
\buildrel {\cong} \over \to
{\rm Hom}_{\Lambda}(\underline{T^{**}}, \underline{H})$.

\vspace{0.2cm}

{\it Proof.} Let $T$ be in mod $\Lambda ^{op}$. Then $T^*$ is
in mod $\Lambda$.
By Proposition 4.2, there are a projective module $P$ and exact sequences
in mod $\Lambda ^{op}$:
$$0 \to M_2^* \to M_1^*(\cong T^{**}\bigoplus P)
\to B_1^* \to 0$$ and
$$0 \to M_3^* \to M_2^* \to B_2^* \to 0$$
with $M_2 \in \Omega ^2({\rm mod}\ \Lambda)$,
$M_3 \in \Omega ^3({\rm mod}\ \Lambda)$, $B_1^*$ projective and
r.pd$_{\Lambda}B_2^* \leq 1$.

Consider the following push-out diagram:

$$\xymatrix{& 0 \ar[d] & 0 \ar[d] & & \\
& M_3^*\ar[d]\ar@{=}[r] & M_3^*\ar[d] & & \\
0 \ar[r] & M_2^* \ar[r]\ar[d] & M_1^* \ar[r]\ar[d]
& B_1^* \ar[r]\ar@{=}[d] & 0\\
0 \ar[r] & B_2^* \ar[r]\ar[d] & Y_1 \ar[r]\ar[d]
& B_1^* \ar[r] & 0\\
& 0 & 0 & & }$$
From the bottom row in the above diagram we know that
r.pd$_{\Lambda}Y_1 \leq 1$. Notice that the first column
$0\to M_3^* \to M_2^* \to B_2^* \to 0$ is dual exact and
the middle row $0\to M_2^* \to M_1^* \to B_1^* \to 0$ splits,
then it is easy to verify that the middle column
$0\to M_3^* \to M_1^* \to Y_1 \to 0$ is dual exact.

We then consider the following push-out diagram:

$$\xymatrix{& 0 \ar[d] & 0 \ar[d] & & \\
& M_4^*\ar[d]\ar@{=}[r] & M_4^*\ar[d] & & \\
0 \ar[r] & M_3^* \ar[r]\ar[d] & M_1^* \ar[r]\ar[d]
& Y_1 \ar[r]\ar@{=}[d] & 0\\
0 \ar[r] & B_3^* \ar[r]\ar[d] & Y_2 \ar[r]\ar[d]
& Y_1 \ar[r] & 0\\
& 0 & 0 & & }$$
where the middle row is the middle column in the former
diagram, the exactness of the first column follows from
Proposition 4.2, $M_4 \in \Omega ^4({\rm mod}\ \Lambda)$ and
r.pd$_{\Lambda}B_3^* \leq 2$.
From the bottom row in the above diagram we know that
r.pd$_{\Lambda}Y_2 \leq 2$. Notice that both the middle row
and the first column in above diagram are dual exact, we then
get the following exact commutative diagram:

$$\xymatrix{M_1^{**} \ar[r]\ar[d] & M_3^{**} \ar[r]\ar[d] & 0 \\
M_4^{**}\ar@{=}[r] & M_4^{**}\ar[d] \\
& 0 & }$$
So $M_1^{**} \to M_4^{**}$ is epic and hence the middle column
in the above diagram is dual exact.

Continuing this process, we finally
get an exact sequence in mod $\Lambda ^{op}$:
$$0 \to X \to T^{**}\bigoplus P (\cong M_1^*) \to Y \to 0
\eqno{(4)}$$
which is dual exact, where
$X=M_k^*$ (where $M_k \in \Omega ^k({\rm mod}\ \Lambda)$) and
r.pd$_{\Lambda}Y \leq k-2$.

Since $M_k \in \Omega ^k({\rm mod}\ \Lambda)$,
Ext$_{\Lambda}^i
(X, \Lambda)=$Ext$_{\Lambda}^i
(M_k^*, \Lambda)=0$ for any $1 \leq i \leq k-2$. So
Ext$_{\Lambda}^i(Y, \Lambda)\cong$
Ext$_{\Lambda}^i
(T^{**}, \Lambda)$ for any $2 \leq i \leq k-2$. On the
other hand, from the fact that Ext$_{\Lambda}^1
(M_k^*, \Lambda)=0$ and the dual exactness of the sequence (4)
we have that Ext$_{\Lambda}^1(Y, \Lambda)\cong$Ext$_{\Lambda}^1
(T^{**}, \Lambda)$ and thus
Ext$_{\Lambda}^i(Y, \Lambda)\cong$Ext$_{\Lambda}^i
(T^{**}, \Lambda)$ for any $1 \leq i \leq k-2$.
So, if $\underline{T^{**}} \to \underline{H}$ is a homomorphism with
r.pd$_{\Lambda}\underline{H} \leq k-2$,
it then follows from [3] Lemma 2.42 that the exact sequence (4)
induces an isomorphism
${\rm Hom}_{\Lambda}(\underline{Y}, \underline{H})
\buildrel {\cong} \over \to
{\rm Hom}_{\Lambda}(\underline{T^{**}}, \underline{H})$.
We are done. $\blacksquare$

\vspace{0.2cm}

Let $\Lambda$ be a commutative noetherian ring and let
$n$ be a non-negative integer and $M$ in
$\Omega ^n({\rm mod}\ \Lambda)$. An Evans-Griffith representation
of $M$ is an exact sequence in mod $\Lambda$:
$$0\to S \to B \to M \to 0$$ where $B$ is an $n$th syzygy of
Ext$_{\Lambda}^{n+1}(D(M), \Lambda)$ and $S$ is in
$\Omega ^{n+2}({\rm mod}\ \Lambda)$
(c.f. [7]). In the case $\Lambda$ is not necessarily commutative
we also call such an exact sequence an Evans-Griffith representation
of $M$.

\vspace{0.2cm}

{\bf Proposition 4.4} {\it Let} $\Lambda$ {\it be a quasi}
$k$-{\it Gorenstein ring. Then, for any} $0\leq d \leq k-1$, {\it
each module in} $\Omega ^d({\rm mod}\ \Lambda)$ {\it has an
Evans-Griffith representation}.

{\it Proof.} Let $M$ be in $\Omega ^d({\rm mod}\ \Lambda)$. By
Theorem there is a projective module $P$ and an exact
sequence in mod $\Lambda$:
$$0\to B \buildrel {\alpha} \over \longrightarrow M\bigoplus P
\buildrel {\beta} \over \longrightarrow C\to 0$$ satisfying the properties
that $B$ is a $d$th syzygy of Ext$_{\Lambda}^{d+1}(D(M), \Lambda)$ and
$C\in \Omega ^{d+1}({\rm mod}\ \Lambda)$.

Let $\gamma$ be the composition: $B \buildrel {\alpha} \over \longrightarrow
M\bigoplus P \buildrel {(1, 0)} \over \longrightarrow M$,
that is, $\gamma =(1, 0)\alpha$. Suppose that $Q
\buildrel {\delta} \over \longrightarrow M \to 0$ is exact in mod
$\Lambda$ with $Q$ projective. Then we have the following commutative
diagram with exact columns and rows:

$$\xymatrix{& 0 \ar[d]
& 0 \ar[d] &  & \\
0 \ar[r] & S\ar[r] \ar[d] & P\bigoplus Q \ar[r]
\ar[d]^{\tiny \left(\begin{array}{cc} 0 & -\delta \\ 1 & 0\\ 0 & 1
\end{array}\right)} & C\ar[r]
\ar@{=}[d] & 0\\
0 \ar[r] & B\bigoplus Q \ar[r]^>>>>>{\binom {\alpha \ 0} {0 \ 1}}
\ar[d]^{(\gamma ,\ \delta)}& M\bigoplus
P\bigoplus Q \ar[r]^>>>>>{(\beta, 0)} \ar[d]^{(1,\ 0,\ \delta)} & C\ar[r] & 0\\
& M \ar[d] \ar@{=}[r] & M \ar[d] &  & \\
& 0 & 0 &  & }
$$ where $S=$Ker$(\gamma, \delta)$.
By the exactness of the first row we have $S\in \Omega
^{d+2}({\rm mod}\ \Lambda)$. Thus the first column
in the above diagram:
$$0 \to S \to B\bigoplus Q \to M \to 0$$ is an
Evans-Griffith representation of $M$.
We are done. $\blacksquare$

\vspace{0.2cm}

The following is an immediate consequence of Proposition 4.4.

\vspace{0.2cm}

{\bf Corollary 4.5} {\it If} $\Lambda$ {\it is a quasi Auslander
ring, then for any non-negative integer} $d$, {\it each module in}
$\Omega ^d({\rm mod}\ \Lambda)$ {\it has an Evans-Griffith
representation}.

\vspace{0.2cm}

If $\Lambda$ is an Auslander ring, that is, $\Lambda$ is
$k$-Gorenstein for all $k$, then the grade condition in Corollary
4.5 is satisfied for $\Lambda$ (see [8] Theorem 3.7). We point out
that a $k$-Gorenstein ring with both left and right self-injective
dimensions being $k$ is an Auslander ring (see [11] Proposition
1).

On the other hand,
by [6] we have that a commutative noetherian ring $\Lambda$ has
finite self-injective dimension if and only if
gradeExt$_{\Lambda}^i(N, \Lambda)\geq i$
for any $N\in$mod $\Lambda$ and $i \geq 1$. So, by Corollary 4.5,
we immediately have the following

\vspace{0.2cm}

{\bf Corollary 4.6} {\it If} $\Lambda$ {\it is a commutative
noetherian with finite self-injective dimension,
then for any non-negative integer} $d$, {\it each module in}
$\Omega ^d({\rm mod}\ \Lambda)$ {\it has an
Evans-Griffith representation}.

\vspace{0.2cm}

Observe that a special instance of Corollary 4.6 was already
considered by Evans and Griffith in [7] Theorem 2.1. They
claimed that if $\Lambda$ is a commutative noetherian local
ring with finite global dimension and contains a field then each
non-free $d$th syzygy of rank $d$ has an
Evans-Griffith representation. Corollary 4.6
generalizes this result to much more general setting.

\vspace{0.5cm}

{\bf Acknowledgements} The author was partially supported by
Specialized Research Fund for the Doctoral Program of Higher
Education. This paper was finished during a visit of the author to
Okayama University from January to June, 2004. The author is
grateful to Prof. Yuji Yoshino for his kind hospitality.

\end{document}